\documentclass[12pt]{article}

\usepackage{xyling}
\usepackage{amsmath}
\usepackage{amsfonts}
\def\aaatree#1#2#3{\Tree { \K{o} \B{d} \\
\K{\mbox{${#1}$}} \B{d} \\
\K{ \mbox{${#2}$}} \B{d} \\
\K{\mbox{${#3}$}}
}}

\def\abctree#1#2#3{\Tree { &  \K{o} \B{dl} \B{d}  \B{dr} \\
\K{\mbox{${#1}$}} & \K{\mbox{${#2}$}} & \K{\mbox{${#3}$}}
}}

\def\abatree#1#2#3{\Tree { &  \K{o} \B{dl} \B{dr} \\
\K{\mbox{${#1}$}} \B{d} & & \K{\mbox{${#2}$}} \\
\K{\mbox{${#3}$}}
}}

\def\abbtree#1#2#3{\Tree { &  \K{o} \B{dl} \B{dr} \\
\K{\mbox{${#1}$}} & & \K{\mbox{${#2}$}} \B{d} \\
       & & \K{\mbox{${#3}$}} }}

\def\axbctree#1#2#3{\Tree { &  \K{o} \B{d} & \\
& \K{\mbox{${#1}$}} \B{dl} \B{dr}  & \\
\K{\mbox{${#2}$}} && \K{\mbox{${#3}$}}
}}

\title{Hopf Algebras of Heap Ordered Trees and Permutations}
\author{Robert L. Grossman \and Richard G. Larson}
\date{November 2, 2007}

\newtheorem{definition}{Definition}
\newtheorem{thm}[definition]{Theorem}
\newtheorem{eg}[definition]{Example}
\newtheorem{rmk}[definition]{Remark}

\def\D{\cal\bf D}

\newcommand{\bminus}{\mathrm{B}_{-}}

\newcommand{\aplus}{\mathrm{A}^{+}}
\newcommand{\hot}[1]{\mathcal{T}_{#1}}
\newcommand{\khotx}[1]{k\mathcal{T}_{#1} }
\newcommand{\khot}{k\mathcal{T}}
\newcommand{\Sp}[1]{\mathfrak{S}_{#1}}
\newcommand{\kS}{k\mathfrak{S} }
\newcommand{\kSx}[1]{k\mathfrak{S}_{#1} }

\newcommand{\rmst}[2]{\mathrm{st}(#1,#2)}
\newcommand{\rmstz}[1]{\mathrm{st}(#1)}

\newcommand{\smsh}{\mathbin\#}
\newcommand{\pf}{\medskip\noindent{\sc Proof: }}
\begin{document}
\maketitle

\begin{abstract}

A standard heap ordered tree with $n+1$ nodes is a finite rooted
tree in which all the nodes except the root are labeled with the
natural numbers between~1 and~$n$, and that satisfies the property
that the labels of the children of a node are all larger than the label
of the node.
Denote the set of standard heap ordered trees with
$n+1$ nodes by $\mathcal{T}_n$.
Let
\[
\khot = \bigoplus_{n\ge0}\khotx{n}.
\]
It is known that there are Hopf algebra structures on
$\khot$.
Let $\Sp{n}$ denote the symmetric group on $n$ symbols.
Let
\[
\kS = \bigoplus_{n\ge0}\kSx{n}.
\]
We give a bialgebra structure on $\kS$,
and show that there is a natural bialgebra isomorphism from $\khot$ to $\kS$.
\end{abstract}

\section{Introduction}

We denote the set of standard heap ordered trees on $n+1$ nodes by
$\mathcal{T}_n$.  Let $\khotx{n}$ denote the vector space over the
field $k$ whose basis is the set of trees in $\mathcal{T}_n$, and let
\[
\khot = \bigoplus_{n\ge0}\khotx{n}.
\]
In \cite{GL:Hopf}, we defined a noncommutative product and a
cocommutative coproduct on $\khot$ that make $\khot$ it a Hopf algebra.
Another product is given in \cite[Ex.~6.2]{GL:Hopf} which gives another Hopf algebra
structure on $\khot$.

Let $\mathfrak{S}_n$ denote the symmetric group on $n$ symbols,
let $\kSx{n}$  denote the vector space over the field $k$ with basis $\mathfrak{S}_n$, and let
\[
\kS = \bigoplus_{n\ge0}\kSx{n}.
\]

In \cite{AS:PermsTrees}, Aguiar and Sottile showed that there is a
filtration on the Malvenuto-Reutenauer Hopf algebra $\kS$ such that
the associated graded dual is isomorphic to $\khot$.
In \cite{MalvenutoReutenauer} Malvenuto and Reutenauer showed
that there is a Hopf algebra structure on $\kS$ that is related
to the Solomon descent algebra \cite{Solomon}.  The Hopf algebra structure
they introduced on $\kS$ is noncommutative, noncocommutative,
self-dual, and graded.

In \cite[Prop.~3.3]{HNT:permtrees} a Hopf algebra structure on $\kS$ is given:
the product is the shifted concatenation of permutations; the coalgebra structure is
given by the fact that the primitive elements are freely generated as a Lie algebra by
the connected permutations (cycles).
It is proved that this Hopf algebra is isomorphic to the Hopf algebra of heap ordered
trees using that $\kS$ is isomorphic to $U(P(\kS))$.
This Hopf algebra corresponds to the Hopf algebra of heap ordered trees with the product
$\odot$ described in \cite[Cor. 6.4]{GL:Hopf}.

In this paper we give another Hopf algebra structure on $\kS$ which closely mirrors the
Hopf algebra structure on heap ordered trees, and which we can easily prove is isomorphic
to that Hopf algebra.

\section{Heap Ordered Trees}\label{HOTrees}

In this section, we define heap ordered trees and standard
heap ordered trees.  In the next section, we define a Hopf algebra
structure on the vector space whose basis is the set of standard heap ordered
trees.

\begin{definition}
  A standard heap ordered tree on $n+1$ nodes is a finite, rooted tree
  in which all nodes except the root are labeled with the numbers
  $\{1$, $2$, $3$,\ldots, $n\}$ so that:
\begin{enumerate}
\item\label{hot:1} each label $i$ occurs precisely once in the tree;
\item\label{hot:2} if a node labeled $i$ has children labeled $j_1$, $\ldots$, $j_k$,
then $i < j_1$,\ldots, $i < j_k$.
\end{enumerate}
\end{definition}

We denote the set of standard heap ordered trees on $n+1$ nodes by $\mathcal{T}_n$.
Let $\khotx{n}$ be
the vector space over the field $k$ whose basis is the set of trees in $\mathcal{T}_n$,
and let
\[
\khot = \bigoplus_{n\ge0}\khotx{n}.
\]

\begin{definition}
A heap ordered tree is a rooted tree in which every node
(including the root) is given a different positive integer label such that
condition~(\ref{hot:2}) is satisfied.
\end{definition}

A heap ordered tree differs from a standard heap ordered tree in that
the root is also labeled, and that the labels can be taken from a
larger set of positive integers.

Our convention for drawing trees is that the root is at the top and
children are drawn so that a child of a node is to the left of another
child if its label is lower.

In the following sections we use the operation that relabels
a heap ordered tree by redefining the labels based on the order they are in and
adding the same positive integer to each.

\begin{definition}
Let $t$ be a heap ordered tree.
Define $\rmst{t}{m}$ (where $m$ is a non-negative integer) as follows:
let $L = (j_1,\ldots,j_k)$ be the ordered list of integers which occur as labels
of nodes of $t$.
The ordered labeled tree $\rmst{t}{m}$ is the same as $t$ except that the
label $j_p$ is replaced by $p+m$.
Note that if the root of $t$ is unlabeled, $\rmst{t}{0}$ is a standard heap ordered tree,
which we denote by $\rmstz{t}$.
\end{definition}

Sometimes it is helpful to think of $\rmst{t}{m}$ as a way
of ``normalizing'' the heap ordered tree $t$.

\section{A Hopf Algebra Structure on Standard Heap Ordered Trees}
In this section, we define a Hopf algebra structure on standard heap ordered trees,
that is, a bialgebra structure which is graded and connected, and so has an antipode.

Suppose that $t_1\in\hot{m}$ and $t_2\in\hot{n}$ are standard heap ordered trees.
Let $s_1$, \ldots, $s_r$ be the children of the root of $\rmst{t_1}{n}$.
We use the operation of deleting the root of a tree to
produce a forest of trees, which we denote
\[
\bminus(\rmst{t_1}{n}) = \{ s_1, \ldots, s_r \}.
\]
Note that the $s_i$ have labeled roots.
If $t_2$ has $n+1$ nodes (counting the root), there are
$(n+1)^r$ ways to attach the $r$ subtrees of $t_1$
which have $s_1$, \ldots, $s_r$ as roots to the tree $t_2$
by making each $s_i$ the child of some node of $t_2$.
We denote the set of copies of $t_2$ with the $s_i$ attached by
\[
\aplus(\{s_1, \ldots, s_r\}, t_2).
\]

The product $t_1t_2$ of the trees $t_1$ and $t_2$ is the
sum of these $(n+1)^r$ trees.  We summarize:

\begin{definition} The product of the two standard heap ordered trees $t_1t_2$ is:
\[
t_1 t_2 = \sum\aplus\left(\bminus(\rmst{t_1}{n}), t_2 \right)\in \khot.
\]
\end{definition}

We define the coalgebra structure of standard heap ordered trees as follows:
\begin{eqnarray*}
\Delta(t) & = &\sum_{X\subseteq\bminus(t)}\rmstz{\aplus(X;e)}\otimes
   \rmstz{\aplus(\bminus(t)\backslash X;e)}  \\
\epsilon(t) & = &
      \begin{cases}
      1 & \text{if $t$ is the tree whose only node is the root,}  \\
      0 & \text{otherwise.}
      \end{cases}
\end{eqnarray*}
Here, if $X\subseteq Y$ are multisets, $Y\backslash X$ denotes the set
theoretic difference.

In \cite{GL:Hopf}, we show that this product and coproduct make $\khot$
into a Hopf algebra.

\section{A Bialgebra Structure for Permutations}
In this section, we define a bialgebra structure on $\kS$.
We begin with some notation.

Let $(\sigma_1 \sigma_2 \cdots \sigma_k)$ denote the cycle in $\Sp{n}$
which sends $\sigma_1$ to $\sigma_2$, $\sigma_2$ to $\sigma_3$, $\ldots$, and
$\sigma_k$ to $\sigma_1$.
Every permutation is a product of disjoint cycles.
If $(\sigma_1 \sigma_2 \cdots \sigma_k)$ is a cycle, then
there is a string naturally associated with the cycle
that we write $\sigma_1 \sigma_2 \cdots \sigma_k$.

Now let $\sigma=(s_1)\cdots(s_r)\in\Sp{m}$ and
$\tau=(t_1)\cdots(t_\ell)\in\Sp{n}$ be two permutations each written
as a product of disjoint cycles.  We denote the corresponding strings
as $s_i=m_{i1}\cdots m_{ip_i}$ and $t_j=n_{j1}\cdots n_{jq_j}$
respectively.  We call the elements of $\{n_{11}, \ldots, n_{\ell
  q_\ell}\}$ {\em attachment points} for the cycles $(s_1)$, \ldots,
$(s_r)$ on the permutation $\tau=(t_1)\cdots(t_\ell)$.  We will also
define $\circ$ to be the $(n+1)^{\mathrm th}$ attachment point.
%In other words, the attachment point $\circ$ is the right of all
%the other attachment points.

Also, if $\sigma$ is a permutation on $\{1,\ldots,k\}$, let
$\rmst{\sigma}{m}$ be the permutation on $\{m+1,\ldots,m+k\}$ that
sends $m+i$ to $m+\sigma(i)$.

The definition of the heap product is simpler if we introduce the standard
order of a permutation, which is defined as follows:
\begin{definition}
We say that a permutation $\sigma\in\Sp{m}$ that is expressed
as a product of cycles $\sigma=(s_1)\cdots(s_r)$ is in standard order if
the cycles $s_i=m_{i1}\cdots m_{ip_i}$ are written so that
\begin{enumerate}
\item $m_{i1} < m_{i2}$, $m_{i1} < m_{i3}$, $\ldots$, $m_{i1} < m_{ip_i}$
\item $m_{11} > m_{21} > m_{31}$, $\ldots$, $m_{i-1,1} > m_{i1}$
\end{enumerate}
\end{definition}

In other words, a product of cycles $\sigma=(s_1)\cdots(s_r)$
is written in standard order if
each cycle $(s_i)=(m_{i1}\cdots m_{ip_i})$ starts with its smallest
entry, and if the cycles $(s_1)$, \ldots, $(s_r)$ are ordered so that their
starting entries are decreasing.  A permutation can always be written in
standard order since disjoint cycles commute, and since a single cycle
is invariant under a cyclic permutation of its string.

We now define the {\em heap product\/} of two permutations.  Given two
permutations $\sigma\in\Sp{m}$ and $\tau\in\Sp{n}$, write them in
standard order $\sigma=(s_1)\cdots(s_r)$ and
$\tau=(t_1)\cdots(t_\ell)$, where the string
$s_i=m_{i1}\cdots m_{ip_i}$ and the string $t_j=n_{j1}\cdots
n_{jq_j}$,

\begin{definition}
We define the heap product $\sigma\smsh\tau$
of $\sigma\in\Sp{m}$ and $\tau\in\Sp{n}$ as follows:
\begin{enumerate}
\item replace $\sigma$ by $\rmst{\sigma}{n}$:
%\item Second, express $\sigma$ and $\tau$ in standard order
%as above.
\item form all terms of the following form:
If $(s_i)$ is one of the cycles in $\sigma$, attach the string $s_i$ to any one
of the $n+1$ attachment points of $\tau$;
if the attachment point is one of $n_{11}$,\ldots, $n_{\ell q_\ell}$, say $n_{jk}$,
place the string $s_i=m_{i1}\cdots m_{ip_i}$ to the right of $n_{jk}$; otherwise
(if the attachment point is $\circ$) we multiply
the term we are constructing by $(s_i)$;
\item The product $\sigma\smsh\tau$ is the sum of all the terms
constructed in this way, taken over all the cycles in $\sigma$
and over all attachment points.
\end{enumerate}
\end{definition}
Note that there are $(n+1)^r$ terms in $\sigma\smsh\tau$.

Some examples will illustrate this.

\begin{eg}
Let $\tau = (n_1 n_2\cdots n_p)\in\Sp{p}$ be a single cycle and
also let $\sigma = (m_1m_2m_3)\in\Sp{3}$ be a single cycle.
We assume that $\{m_1,m_2,m_3\}=\{p+1,p+2,p+3\}$.
We compute the
product $\tau\smsh\sigma$  as follows:
\begin{eqnarray*}
\sigma\smsh\tau & = &
(n_1m_1m_2m_3n_2\cdots n_p) + (n_1n_2m_1m_2m_3\cdots n_p)+\cdots \\
  &  & {}+(n_1n_2\cdots n_pm_1m_2m_3)+(n_1n_2\cdots n_p)(m_1m_2m_3),
\end{eqnarray*}
giving $p+1$ terms.
\end{eg}

\begin{eg}
Let $\sigma=(m_1)(m_2)(m_3)\in\Sp{3}$ be the product of three 1-cycles,
and let $\tau=(n_1n_2n_3)\in\Sp{3}$ be a 3-cycle.
We assume that $\{m_1,m_2,m_3\}=\{4,5,6\}$.
Then
\begin{eqnarray*}
\sigma\smsh\tau & = & (n_1m_1m_2m_3n_2n_3)+(n_1m_1m_2n_2m_3n_3)+(n_1m_1m_2n_2n_3m_3) \\
&  & {}+(n_1m_1m_3n_2m_2n_3)+(n_1m_1n_2m_2m_3n_3)+(n_1m_1n_2m_2n_3m_3) \\
&  & {}+\cdots  \\
&  & {}+(n_1m_3n_2n_3)(m_1)(m_2)+(n_1n_2m_3n_3)(m_1)(m_2) \\
&  & {}+(n_1n_2n_3m_3)(m_1)(m_2)+(n_1n_2n_3)(m_1)(m_2)(m_3)
\end{eqnarray*}
giving $4^3=64$ terms.
\end{eg}

We now describe the coalgebra structure of $\kS$.

We define a function $\rmstz{\pi}$ from permutations to permutations as follows:
let $\pi=(s_1)\cdots(s_p)\in\Sp{n}$ and
let let $L=\{\ell_1,\ldots,\ell_k\}$ be the labels (in order) which occur in the $s_i$.
(If $\pi$ fixes $i$, we include a 1-cycle $(i)$ as a factor in $\pi$.)
The permutation $\rmstz{\pi}$ is the permutation in $\Sp{k}$ gotten by replacing $\ell_j$
with $j$ in $\pi$.
For example, if $\pi=(1 3)(4)(5 7)\in\Sp{7}$, then $\rmstz{\pi}\in\Sp{5}$ equals
$(1 2)(3)(4 5)$.

The coalgebra structure of $\kS$ is defined as follows.
let $\pi=(s_1)\cdots(s_k)\in\Sp{n}$, and let $C=\{(s_1),\ldots,(s_k)\}$.
If $X\subseteq C$ let
$\rho(X) = \rmstz{\prod_{(s_i)\in C}(s_i)}$.
Note that if $\rho(X)\in\Sp{k}$, then $\rho(C\backslash X)\in\Sp{n-k}$.
Define
\begin{eqnarray*}
\Delta(\pi) & = &\sum_{X\subseteq C}\rho(X)\otimes
   \rho(C\backslash X)  \\
\epsilon(\pi) & = &
      \begin{cases}
      1 & \text{if $\pi$ is the identity permutation in $\Sp{0}$,}  \\
      0 & \text{otherwise.}
      \end{cases}
\end{eqnarray*}

In Section~\ref{sectionHOT2P}, we show that there is a bialgebra
isomorphism between the Hopf algebra of standard heap ordered trees and
the bialgebra of permutations.

\section{From Standard Heap Ordered Trees to Permutations}\label{sectionHOT2P}

In this section, we define a map $\varphi$ from
standard heap ordered trees $\hot{n}$ to permutations $\Sp{n}$ and
show that this gives a bialgebra homomorphism from $\khot$ to $\kS$.
We show that $\varphi$ is an isomorphism by giving its inverse.

\begin{definition}
We define the map $\alpha$ from heap ordered trees to strings recursively.
Let $t$ be a heap ordered tree:
if the root of $t$ has label $i$ and children $t_1$, \ldots,
$t_k$ (read from left to right), then
$\alpha(t)$ is the string $i\alpha(t_k)\cdots\alpha(t_1)$.
Note that if $t$ (whose root is labeled with $i$) has no children,
then $\alpha(t)$ is the string $i$.
\end{definition}

It is important to note that this definition relies on the convention
defined in Section~\ref{HOTrees}.  With this convention, children of a node
are arranged from left to right in {\em increasing} order.

Let $t$ be a standard heap ordered tree, and let $t_1$, \ldots, $t_k$ be the
heap ordered trees which are the children of the root of $t$.
Then $\varphi(t)$ is the permutation $(\alpha(t_1))\cdots(\alpha(t_k))$.

\begin{definition}
We recursively define a map $\beta$ from strings to heap ordered trees.
We define $\beta$ on strings whose elements are either numbers or heap ordered trees.
In comparisons which involve heap ordered trees, we use the label of the root of
the tree in the comparison.

Let $s$ be a string $n_1\cdots n_k$ with $n_1<n_2$, \ldots, $n_1<n_k$.
A {\em valid\/} substring is a substring $n_i\cdots n_j$ with $n_i<n_{i+1}$, \ldots,
$n_i<n_j$ and either $n_i>n_{j+1}$ or $j=k$ (that is, $n_j$ is the last entry in $s$).
A string always has a valid substring (which might be the whole string).

The map $\beta$ is defined as follows:
let $n_i\cdots n_j$ be a valid substring.
Replace this substring with the heap ordered tree with root labeled with $n_i$,
and with children of the root $n_k$ ($i<k\le j$) either labeled with $n_k$ (if $n_k$ is
a number) or the tree itself (if $n_k$ is a heap ordered tree).
\end{definition}

Now we can define $\varphi^{-1}$:
Let $\pi=(s_1)\cdots(s_k)$ be a permutation in $\Sp{n}$.
If some number $i$ in $\{1,\ldots,n\}$ does not occur in any of the strings $s_1$, \dots,
$s_k$, replace $\pi$ by $\pi(i)$ so that we may assume that every number
in $\{1,\ldots,n\}$ occurs explicitly in $\pi$.
We can assume that $\pi$ is in standard order.

We construct a standard heap ordered tree $t=\varphi^{-1}(\pi)$ as follows:
first construct the root and let $\beta(s_1)$, \ldots, $\beta(s_k)$ be the children
of that unlabeled root.

\begin{thm}
The map $\varphi$ is a bialgebra morphism.
\end{thm}

\pf
In forming the product in $\khot$ we attach the children of the root of the first
multiplicand to nodes of the second multiplicand or to its root.
In forming the product in $\kS$ we attach cycles of the first
multiplicand to attachment points of the second multiplicand.
These are essentially the same operation.

In the construction of the coproduct of $\khot$, for $t\in\hot{n}$ we use each subset of
the set of children of the root on $t$ to construct a new tree.
In the construction of the coproduct on $\kS$ for $\pi\in\Sp{n}$ we use each subset of
the set of cycles of $\pi$ to construct a new permutation.
These are essentially the same operation.

\medskip

\begin{rmk}
Since we already know that $\khot$ is a Hopf algebra with associative multiplication
and that $\kS\cong\khot$, this gives a proof that $\kS$ is a Hopf algebra.
\end{rmk}

Here are some examples of $\varphi$:

\begin{eqnarray*}
\mbox{\abctree{1}{2}{3}}  & \stackrel{\varphi}{\longrightarrow} & (1)(2)(3) \\
\mbox{\aaatree{1}{2}{3}}  & \stackrel{\varphi}{\longrightarrow} & (123) \\
\end{eqnarray*}
\begin{eqnarray*}
\mbox{\abatree{1}{2}{3}}  & \stackrel{\varphi}{\longrightarrow} & (13)(2) \\
\mbox{\abbtree{1}{2}{3} } & \stackrel{\varphi}{\longrightarrow} & (1)(23) \\
\mbox{\axbctree{1}{2}{3}} & \stackrel{\varphi}{\longrightarrow} & (132)
\end{eqnarray*}

\end{document}